\documentclass[10pt,leqno]{amsart}
\usepackage{graphicx}
\usepackage{titlesec}

\usepackage{indentfirst,csquotes}

\usepackage{amssymb,amsthm,amsmath}
\usepackage{xcolor,paralist, hyperref, titlesec, fancyhdr, etoolbox}
\usepackage{titlesec}
\newtheorem{theorem}{Theorem} 

\newtheorem{lem}[theorem]{Lemma} 

\newtheorem{corollary}[theorem]{Corollary} 

\theoremstyle{definition} 
\theoremstyle{remark} 
\newtheorem{remark}[theorem]{Remark}

\titleformat{\section}[display]{\normalfont\huge\bfseries\centering}{\centering\chaptertitlename\thechapter}{10pt}{\Large}
\titlespacing*{\section}{2pt}{2ex}{2ex}

\hypersetup{ colorlinks=true, linkcolor=black, filecolor=black, urlcolor=black }

\usepackage{lipsum}
\begin{document}
\title{Generalization of Some Well-Known Polynomial Inequalities for the Modified Smirnov Operator} 
\author[Deepak Kumar, Dinesh Tripathi, Sunil Hans]{
  Deepak Kumar\textsuperscript{1}, Dinesh Tripathi\textsuperscript{2}, Sunil Hans\textsuperscript{1}
}
\date{\today}
\address{Address}
\address{Author (1): Department of Applied Mathematics, Amity University,\\ Noida-201313, India.}
\email{deepak.kumar26@s.amity.edu}
\email{sunil.hans82@yahoo.com}
\address{Author (2): Department of Science - Mathematics, School of Sciences, Manav Rachana University, Faridabad-121004, India.}
\email{dinesh@mru.edu.in}
\date{}
\let\thefootnote\relax
\footnotetext{%
}
\begin{abstract}
  Let $P(z)$ be a polynomial of degree $n$.  In this paper, we consider the modified Smirnov operator, which carries a polynomial $P(z)$ into $\tilde{\mathbb{S}}_a[P](z)=(1+az)P'(z)-naP(z),$ where $a$ is an arbitrary number in $\overline{\mathbb{D}}$. We estimate minimum and maximum moduli of modified Smirnov operator of $P(z)$ on the unit circle with restricted zeros and thereby obtain a generalization of some results of Dewan and Hans \cite{dewan2013some}. This study includes compact generalization of some well-known polynomial inequalities.
\end{abstract} 
\maketitle
\begin{center}
\section*{INTRODUCTION}
\end{center}
 Let $P(z)$ be a polynomial of degree $n$ and $P'(z)$ be the derivative of polynomial $P(z)$. Let $\mathbb{D}$ be the open unit disk $\{z\in\mathbb{C}; |z|<1\} $, so that $\overline{\mathbb{D}}$ is it's closure and $B(\mathbb{D) }$ denotes its boundary, then
     \begin{align}{\label{1.1}}
         \max_{z\in B(\mathbb{D})}|P'(z)| \leq n \max_{z\in B(\mathbb{D})}|P(z)|,
     \end{align}
     and
     \begin{align}{\label{1.2}}
         \max_{z\in B(\mathbb{D})}|P(Rz)| \leq R^n \max_{z\in B(\mathbb{D})}|P(z)|.
     \end{align}
    Inequality (\ref{1.1}) can be obtained by a direct result of a theorem of S. Berstein \cite{bernstein1912ordre} for the derivative of a polynomial. A straightforward inference from the maximum modulus principle \cite{polya1925aufgaben} yields the inequality (\ref{1.2}). In both inequalities, the equality holds for $P(z)=\lambda z^n, \lambda\neq 0$. 
    Aziz and Dawood \cite{AZIZ1988306} have shown that if $P(z)$ has all its zeros in $\overline{\mathbb{D}}$, then
    \begin{align}{\label{3}}
        \min_{z \in B(\mathbb{D})}|P'(z)|\geq n\min_{z \in B(\mathbb{D})}|P(z)|,
    \end{align}
    and 
    \begin{align}{\label{4}}
        \min_{z \in B(\mathbb{D})}|P(Rz)| \geq R^n \min_{z \in B(\mathbb{D})}|P(z)|.
    \end{align}
    Inequalities (\ref{1.1}) and (\ref{1.2}) can be improved if the zeros are restricted. Erd\"os conjectured and Lax \cite{lax1944proof} proved that if $P(z)$ has no zeros in $\mathbb{D}$, then
\begin{align}{\label{1.3}}
        \max_{z\in B(\mathbb{D})}|P'(z)| \leq \frac{n}{2}\max_{z\in B(\mathbb{D})}|P(z)|,
    \end{align}
    and for $R \geq 1$,
    \begin{align}{\label{1.4}}
        \max_{z\in B(\mathbb{D})}|P(Rz)| \leq \frac{R^n+1}{2}\max_{z\in B(\mathbb{D})}|P(z)|.
    \end{align}
   Ankeny and Rivlin \cite{ankeny1955theorem}  used (\ref{1.3}) to prove the above inequality (\ref{1.4}).
   Inequalities (\ref{1.3}) and (\ref{1.4}) were further improved in \cite{AZIZ1988306}, where under the same hypothesis, it was shown that 
   \begin{align}{\label{5}}
       \max_{z \in B(\mathbb{D})}|P'(z)| \leq \frac{n}{2}\biggl\{\max_{z \in B(\mathbb{D})}|P(z)|-\min_{z \in B(\mathbb{D})}|P(z)|\biggr\},
   \end{align}
   and 
   \begin{align}{\label{6}}
       \max_{z \in B(\mathbb{D})}|P(Rz)| \leq \left( \frac{R^n+1}{2}\right)\max_{z \in B(\mathbb{D})}|P(z)|- \left(\frac{R^n-1}{2}\right)\min_{z \in B(\mathbb{D})}|P(z)|.
   \end{align}
   The equality (\ref{1.3}), (\ref{1.4}), (\ref{5}) and (\ref{6}) holds for the polynomial $P(z)=\alpha z^n+\beta$, where $|\alpha|=|\beta|.$ Some other generalization of inequalities (\ref{5}) and (\ref{6}) has been mentioned in \cite{dt1}, \cite{dt2} and \cite{dt3}.\\
   Aziz and Rather \cite{aziz1999inequality} generalizes the inequality  (\ref{1.1}) and (\ref{1.2}) and proved that if $P(z)$ is a polynomial of degree $n$, then for every real and complex number $\alpha$ with $|\alpha|\leq 1$ and $R\geq 1$,
   \begin{align}{\label{1.5}}
       |P(Rz)-\alpha P(z)|\leq |R^n-\alpha||z|^n\max_{z \in B(\mathbb{D})}|P(z)| \ \text{for} \ z \in \mathbb{C\backslash D}.
   \end{align}
   The above result is best possible and equality in (\ref{1.5}) holds for $P(z)=\lambda z^n, \lambda\neq 0$. Inequality (\ref{1.1}) follows directly from (\ref{1.5}) by dividing both sides by $R-1$ and taking limit $R \to 1$ with $\alpha=1$. Inequality (\ref{1.5}) reduces to (\ref{1.2}) for $\alpha =0$.
      As an improvement of inequality (\ref{1.5}), the authors \cite{aziz1999inequality} have also shown that if $P(z)\neq 0$ in $\mathbb{D}$, then for every real and complex number $\alpha$ with $|\alpha|\leq 1$ and $R\geq 1$
   \begin{align}{\label{1.6}}
        |P(Rz)-\alpha P(z)|\leq \frac{1}{2}\{|R^n-\alpha||z|^n+|1-\alpha|\}\max_{z \in B(\mathbb{D})}|P(z)|,
   \end{align}
   for $ z \in \mathbb{C\backslash D}$.
   The result is sharp and equality in (\ref{1.6}) holds for $P(z)=z^n + 1$. After dividing both sides of (\ref{1.6}) by $(R-1)$ and taking limit $R \to 1$ with $\alpha=1$, inequality (\ref{1.3}) is obtained. Furthermore, inequality (\ref{1.6}) reduces to (\ref{1.4}) for $\alpha =0$. \\
 Aziz and rather \cite{AZIZ1988306} generalized the inequalities (\ref{3}) and (\ref{4}) and demonstrated that if $P(z) \neq 0$ in $\mathbb{C\backslash D}$, then for every real or complex number $\alpha$ with $|\alpha|\leq 1$ and $R\geq 1$
   \begin{align}{\label{7}}
       \min_{z \in B(\mathbb{D})}|P(Rz)-\alpha P(z)| \geq |R^n-\alpha|\min_{z \in B(\mathbb{D})}|P(z)|.
   \end{align}
   The result is sharp and equality holds for $P(z)=me^{i\gamma}z^n, m>0 $. As an improvement in inequality (\ref{1.6}), Aziz and Rather \cite{AZIZ1988306} have also proved that if $P(z)\neq 0 $ in $\mathbb{D}$, then for every real or complex number $\alpha$ with $|\alpha|\leq 1$ and $R\geq 1$
\begin{align}{\label{8}}
\nonumber|P(Rz)-\alpha P(z)| \leq & \frac{1}{2}\biggl[\{|R^n-\alpha| +|1-\alpha|\}\max_{z \in B(\mathbb{D})}|P(z)| \\ &- \{|R^n-\alpha|-|1-\alpha|\}\min_{z \in B(\mathbb{D})}|P(z)|\biggr].
   \end{align}
   Inequality (\ref{8}) is sharp and equality holds for $P(z)=\beta+\gamma z^n$, where $|\beta|=|\gamma|=\frac{1}{2}$.\\
   In $1930$, S. Bernstein \cite{bernstein1930limitation} also proved the following result:
   \begin{theorem}{\label{thmA}}
       Let $F(z)$ be a polynomial of degree $n$, having all its zeros in $\overline{\mathbb{D}}$ and $P(z)$ be a polynomial of degree not exceeding that of $F(z)$. If $|P(z)|\leq |F(z)|$  on $B(\mathbb{D})$, then 
       \begin{align*}
          |P'(z)| \leq |F'(z)| \ for \ z\in \mathbb{C} \backslash \mathbb{D}.
       \end{align*} 
       The equality holds only if $P(z)=e^{i\gamma}F(z), \gamma \in \mathbb{R}$.
   \end{theorem}
   For $z\in \mathbb{C}\backslash \mathbb{D}$, denoting $\Omega _{|z|}$ the image of the disk $\{t\in \mathbb{C}; |t|<|z| \}$ under the mapping $\phi(t)= \frac{t}{t+1}$, Smirnov \cite{smirnov1964constructive} as a generalization of Theorem \ref{thmA} proved the following:
   \begin{theorem}{\label{B}}
      Let $F(z)$ and $P(z)$ be the polynomials possessing condition as in Theorem \ref{thmA}. Then  for $z\in \mathbb{C\backslash D}$
      \begin{align}{\label{1.7}}
          |\mathbb{S}_\alpha[P](z)| \leq |\mathbb{S}_\alpha[F](z)|,
            \end{align}
          for all $\alpha \in \overline{\Omega}_{|z|}$, with $\mathbb{S}_\alpha[P](z)= z P'(z)-n\alpha P(z)$, where $\alpha$ is a constant.  
   \end{theorem}
   For $\alpha \in \overline{\Omega}_{|z|}$ in the inequality (\ref{1.7}), the equality is holds at a point $z \in \mathbb{C\backslash D}$, only if $P(z) = e^{i \gamma}F(z), \gamma \in \mathbb{R}$. We note that for fixed $z \in \mathbb{C \backslash D}$, the inequality (\ref{1.7}) can be replaced by (see for reference \cite{ganenkova2019variations})
   \[\biggl|zP'(z)-n\frac{az}{1+az}P(z)\biggr| \leq \biggl|zF'(z)-n\frac{az}{1+az}F(z)\biggr|,\]
   where $a$ is an arbitrary number in $\overline{\mathbb{D}}$.\\
   Equivalently for $z \in \mathbb{C\backslash D}$
   \[|\Tilde{\mathbb{S}}_a[P](z)| \leq |\Tilde{\mathbb{S}}_a[F](z)|,\]
   where $\Tilde{\mathbb{S}}_a[P](z)=(1+az)P'(z)-naP(z)$ is known as the modified Smirnov operator. The modified Smirnov operator $\Tilde{\mathbb{S}}_a$ is preferred over Smirnov operator $\mathbb{S}_\alpha$, because the parameter $a$ of $\Tilde{\mathbb{S}}_a$ does not depend on $z$ unlike the parameter $\alpha$ of $\mathbb{S}_\alpha$.\\
   Shah and Fatima \cite{shah2022bernstein} used modified Smirnov operator to generalize inequalities (\ref{1.1}), (\ref{1.2}), (\ref{1.3}) and (\ref{1.4}) and proved that if $P(z)$ is a polynomial of degree $n$, such that $|P(z)| \leq M$ for $z\in B(\mathbb{D})$, then 
 for $ z \in \mathbb{C\backslash D}$
   \begin{align}{\label{1.8}}
      | \tilde{\mathbb{S}}_a[P](z)|\leq M|\tilde{\mathbb{S}}_a[z^n]|,
       \end{align}
       equivalently
       \begin{align}{\label{1.9}}
           |(1+az)P'(z)-naP(z)| \leq Mn|z|^{n-1},
       \end{align}
and if $P(z)\neq 0 $ in $\mathbb{D}$, then for $ z \in\mathbb{C\backslash D}$
   \begin{align}{\label{1.10}}
    |\tilde{\mathbb{S}}_a[P](z)|\leq \frac{1}{2}\{|\tilde{\mathbb{S}}_a[z^n]|+n|a|\}\max_{z\in B(\mathbb{D})}|P(z)|,
   \end{align}
   equivalently
   \begin{align}\label{1.11}
       |(1+az)P'(z)-naP(z)| \leq \frac{1}{2}\{n|z|^{n-1}+n|a|\}\max_{z \in B(\mathbb{D})}|P(z)|.
   \end{align}
   By setting $a=0$, inequalities $(\ref{1.9})$ and (\ref{1.11}) reduce to inequalities (\ref{1.1}) and (\ref{1.3}) respectively. Likewise, by setting $ a = -\frac{1}{z}$ with $z = Re^{i\theta}, R \geq 1$, inequalities (\ref{1.2}) and (\ref{1.4}) can be obtained from inequalities (\ref{1.9}) and (\ref{1.11}), respectively.\\
   The authors \cite{shah2022bernstein} applied a modified Smirnov operator to generalize the inequalities (\ref{3}), (\ref{4}), (\ref{5}) and (\ref{6}) and demonstrated that if $P(z)$ is a polynomial of degree $n$ with $P(z)\neq 0$ in $\mathbb{C\backslash \overline{D}}$, then 
   \begin{align}{\label{9}}
       |\tilde{\mathbb{S}}_a[P](z)| \geq |\tilde{\mathbb{S}}_a[z^n]| \min_{z \in B(\mathbb{D})}|P(z)|,
   \end{align} 
   and if $P(z)\neq 0$ for $ z \in \mathbb{D}$, then for $z \in \mathbb{C\backslash D}$
   \begin{align}{\label{10}}
        |\tilde{\mathbb{S}}_a[P](z)| \leq \frac{1}{2}\{|\tilde{\mathbb{S}}_a[z^n]| +n|a|\} \max_{z \in B(\mathbb{D})}|P(z)|- \frac{1}{2}\{|\tilde{\mathbb{S}}_a[z^n]| -n|a|\} \min_{z \in B(\mathbb{D})}|P(z)| .
   \end{align}
   The above result is best possible and equality holds for the polynomials having all zeros on unit disk. \\
   Wani and Liman \cite{wani2024bernstein} have generalized inequality (\ref{1.5}) concerning the modified Smirnov operator and proved that if $P(z)$ is a polynomial of degree $n$, then for every real and complex number $\alpha$ with $|\alpha| \leq 1$ and $R \geq 1$ 
\begin{align}{\label{1.12}}
       |\tilde{\mathbb{S}}_a[P](Rz)-\alpha \tilde{\mathbb{S}}_a[P](z)| \leq |R^n-\alpha||\tilde{\mathbb{S}}_a[z^n]|\max_{z\in B(\mathbb{D})}|P(z)|,
   \end{align}
    for $z \in \mathbb{C\backslash D}$. The result is sharp and equality holds in (\ref{1.12}) for $\lambda z^n, \lambda \neq 0$.
   Further, as a generalization of inequality (\ref{1.6}), the authors \cite{wani2024bernstein} have also shown that if $P(z) \neq 0$ in $\mathbb{D}$, then for every real and complex number $\alpha$ with $|\alpha|\leq 1$ and $R \geq 1$
   \begin{align}{\label{1.13}}
       |\tilde{\mathbb{S}}_a[P](Rz)-\alpha\tilde{\mathbb{S}}_a[P](z)| \leq \biggl\{\frac{|R^n-\alpha||\tilde{\mathbb{S}}_a[z^n]|+n|1-\alpha||a|}{2}\biggr\}\max_{z \in B(\mathbb{D})}|P(z)|, 
   \end{align}
    for $\ z\in \mathbb{C\backslash D}$. The result is the best possible and the equality holds for $P(z)= \lambda z^n, \lambda \neq 0$.
 As a compact generalization of the inequalities (\ref{3}), (\ref{4}) and (\ref{7}), Dewan and Hans \cite{dewan2013some} proved the following theorem:
   \begin{theorem}{\label{thmC}}
       If P(z) is a polynomial of degree n, having all its zeros in $\mathbb{D}$, then for every real and complex number $ \alpha,\beta$ with $|\alpha| \leq 1, |\beta| \leq 1 \  and \ R \geq 1$,
       \begin{align}{\label{1.14}}
          \nonumber &\min_{z \in B(\mathbb{D})}\bigg|P(Rz)-\alpha P(z)+\beta\biggl\{\biggl(\frac{R+1}{2}\biggr)^n-|\alpha|\biggr\}P(z)\bigg|\\& \leq \bigg|R^n-\alpha+\beta\biggl\{\biggl(\frac{R+1}{2}\biggr)^n-|\alpha|\biggr\}\bigg|\min_{z \in B(\mathbb{D})}|P(z)|.
       \end{align}
        The result is sharp and the equality in \textnormal{(\ref{1.14})} holds for $P(z)=me^{i\gamma} z^n, m > 0$.
   \end{theorem}
    As an improvement of the above result, Dewan and Hans \cite{dewan2013some} proved the following theorem for the class of polynomials having no zeros in the unit disk and obtain generalization of the inequalities (\ref{5}), (\ref{6}), and (\ref{8}).
    \begin{theorem}{\label{thmD}}
        If $P(z)$ is a polynomial of degree $n$, which does not vanish in $\mathbb{D}$, then for every real and complex number $\alpha, \beta$ with $|\alpha| \leq 1, |\beta| \leq 1$ and $R\geq 1$,
        \begin{align}{\label{1.15}}
            \nonumber &\bigg|P(Rz)-\alpha P(z)+\beta\biggl\{\biggl(\frac{R+1}{2}\biggr)^n-|\alpha|\biggr\}P(z)\bigg|\\& \nonumber\leq \frac{1}{2}\biggl[\biggl\{\bigg|R^n-\alpha+\beta\biggl\{\biggl(\frac{R+1}{2}\biggr)^n-|\alpha|\biggr\}\bigg| + \\& \nonumber\bigg|1-\alpha + \beta\biggl\{\bigg(\frac{R+1}{2}\bigg)^n-|\alpha|\biggr\}\bigg|\biggr\}\max_{z \in B(\mathbb{D})}|P(z)| - \\& \nonumber \biggl\{\bigg|R^n-\alpha+\beta\biggl\{\biggl(\frac{R+1}{2}\biggr)^n-|\alpha|\biggr\}\bigg| - \\&\bigg|1-\alpha + \beta\biggl\{\bigg(\frac{R+1}{2}\bigg)^n-|\alpha|\biggr\}\bigg|\biggr\}\min_{z \in B(\mathbb{D})}|P(z)|\bigg].
        \end{align}
        The result is best possible and the equality in\textnormal{ (\ref{1.15})} is holds for $P(z)= \gamma z^n+\delta$ where $|\gamma|=|\delta|=\frac{1}{2}$.
    \end{theorem}
    In this paper, we prove the following results, which generalize Theorem \ref{thmC} and Theorem \ref{thmD} for the modified Smirnov operator.
\begin{center}
  \section*{MAIN RESULTS}	
\end{center}
\begin{theorem}{\label{thm1}}
         If $P(z)$ is a polynomial of degree $n$, having all its zeros in $\overline{\mathbb{D}}$, then for every real and complex number $\alpha, \beta$ with $|\alpha|\leq 1, |\beta|\leq 1$ and $R\geq 1$,
         \begin{align}{\label{2.1}}
            \nonumber &\biggl| \Tilde{\mathbb{S}}_a[P](Rz)-\alpha \Tilde{\mathbb{S}}_a[P](z)+\beta \biggl\{\biggl(\frac{R+1}{2}\biggr)^n-|\alpha| \biggr\}\Tilde{\mathbb{S}}_a[P](z)\biggr|\\& \geq \biggl|R^n-\alpha +\beta\biggl\{\biggl(\frac{R+1}{2}\biggr)^n -|\alpha|\biggr\}\biggr||\tilde{\mathbb{S}}_a[z^n]|\min_{z\in B(\mathbb{D})}|P(z)|,
         \end{align}
        for $\  z \in \mathbb{C \backslash D}$. The result is sharp and equality in \textnormal{(\ref{2.1})} holds for $P(z)=\lambda z^n, \lambda\neq 0.$
     \end{theorem}
     \begin{remark}
        If in inequality (\ref{2.1}), we take $a=0$, we get the following result
        \begin{align*}
        \biggl|& RP'(Rz)-\alpha P'(z)+ \beta\biggl\{\biggl(\frac{R+1}{2}\biggr)^n -|\alpha|\biggr\}P'(z)\bigg|\\&  \geq \bigg|R^n-\alpha+\beta \biggl\{ \biggl( \frac{R+1}{2}\biggr)^n-|\alpha|\biggr\}\bigg|n|z|^{n-1}\min_{z \in B(\mathbb{D})}|P(z)|, \ for \ z\in \mathbb{C\backslash D}
         \end{align*}
        Inequality (\ref{1.14}) is a special case of inequality (\ref{2.1}) for $ a=-\frac{1}{z}$.
     \end{remark}
     \begin{remark}
     If we consider $\alpha$ and $\beta$ to be zero in inequality (\ref{2.1}), then we get
      \begin{align*}
    |\tilde{\mathbb{S}}_a[P](Rz)|\geq R^n |\tilde{\mathbb{S}}_a[z^n]|\min_{z \in B(\mathbb{D})}|P(z)|, \ for \ z\in \mathbb{C\backslash D}.
      \end{align*}
      By substituting, $a = -\frac{1}{z}$ and $R= 1$ in the above inequality, we obtain inequalities (\ref{4}) and a special case of (\ref{9}), respectively.
     \end{remark}
     If we choose $\alpha =1$ in inequality (\ref{2.1}) and on dividing both sides by $R-1$ and taking $ R \to 1$, we get
     \begin{corollary}{\label{cor2.1}}
         If $P(z)$ is a polynomial of degree $n$, having all its zeros in $\overline{\mathbb{D}}$, then for every real and complex number $ \beta$ with $|\beta|\leq 1$ and $R\geq 1$,
         \begin{align}{\label{2.2}}
     \bigg|z\tilde{\mathbb{S}}_a[P'](z)+\frac{n}{2}  \beta  
  \tilde{\mathbb{S}}_a[P](z) +P'(z)\biggr|\geq  n\bigg|1+\frac{\beta}{2}\bigg||\tilde{\mathbb{S}}_a[z^n]|\min_{z \in B(\mathbb{D})}|P(z)|, 
         \end{align}
         for $z\in \mathbb{C\backslash D}$. The equality holds in the above inequality \textnormal{(\ref{2.2})} for $  P(z)=\lambda z^n, \lambda \neq 0$.       
     \end{corollary}
      \noindent For the $ \alpha =0$ inequality (\ref{2.1}), yield the following result, 
     \begin{corollary}
     If $P(z)$ is a polynomial of degree $n$, having all its zeros in $\overline{\mathbb{D}}$,  then for every real and complex number $ \beta$ with $|\beta|\leq 1$ and $R\geq 1$,
     \begin{align*}
         \nonumber \biggl| \Tilde{\mathbb{S}}_a[P](Rz)+\beta \biggl(\frac{R+1}{2}\biggr)^n\Tilde{\mathbb{S}}_a[P](z)\biggr| \geq  \biggl|R^n +\beta\biggl(&\frac{R+1}{2}\biggr)^n \biggr| |\tilde{\mathbb{S}}_a[z^n]| \min_{z\in B(\mathbb{D})}|P(z)|,
         \end{align*}
          for  $z\in \mathbb{C \backslash D}$. The above result is sharp and equality holds for $ P(z)=\lambda z^n, \lambda \neq 0.$
     \end{corollary} 
    \noindent If we take $\beta=0$ in (\ref{2.1}), we get
     \begin{corollary}
     If $P(z)$ is a polynomial of degree $n$, having all its zeros in $\overline{\mathbb{D}}$, then for every real and complex number $\alpha$ with $|\alpha| \leq 1$ and $R\geq 1$
     \begin{align*}
         |\tilde{\mathbb{S}}_a[P](Rz)-\alpha \tilde{\mathbb{S}}_a[P](z)| \geq |R^n-\alpha| |\tilde{\mathbb{S}}_a[z^n]|\min_{z \in B(\mathbb{D})}|P(z)|.
     \end{align*}
     This result is sharp and equality holds for $P(z)= \lambda z^n, \lambda \neq 0$.
      \end{corollary}
     \begin{theorem}{\label{thm2}}
         If $P(z)$ is a polynomial of degree $n$, which doesn't vanish in $\mathbb{D}$, then for every real and  complex number $\alpha, \beta$ with $|\alpha|\leq1,|\beta|\leq 1$ and $R\geq 1$,
         \begin{align}{\label{2.3}}
           \nonumber  \biggl|&\Tilde{\mathbb{S}}_a[P](Rz)-\alpha \Tilde{\mathbb{S}}_a[P](z)+\beta \biggl\{\biggl(\frac{R+1}{2}\biggr)^n -|\alpha|\biggr\}\Tilde{\mathbb{S}}_a[P](z)\biggr|\\& \nonumber\leq \frac{1}{2}\biggl[ \biggl\{\biggl|R^n-\alpha+\beta\biggl\{\bigg(\frac{R+1}{2}\bigg)^n-|\alpha|\biggr\}\biggr||\Tilde{\mathbb{S}}_a[z^n]|+\\& \nonumber\qquad\biggl|1-\alpha + \beta\bigg\{\biggl(\frac{R+1}{2}\biggr)^n-|\alpha|\bigg\}\biggr|n|a|\biggr\}\max_{z \in B(\mathbb{D})}|P(z)| - \\ &\nonumber \quad \biggl\{\bigg|R^n-\alpha+\beta\biggl\{\biggl(\frac{R+1}{2}\biggr)^n-|\alpha|\biggr\}\bigg||\tilde{\mathbb{S}}_a[z^n]|- \\ & \qquad\bigg|1-\alpha+\beta \biggl\{\biggl(\frac{R+1}{2}\biggr)^n-|\alpha|\biggr\}\bigg|n|a|\biggr\}\min_{z \in B(\mathbb{D})}|P(z)|\biggr],
         \end{align}
        for $z\in \mathbb{C\backslash D}$. The equality in \textnormal{(\ref{2.3})} holds for $P(z)=z^n+1$.
     \end{theorem}
   \noindent  If we take $\beta = 0$ in (\ref{2.3}), we get 
     \begin{corollary}
         If $P(z)$ is a polynomial of degree $n$, which doesn't vanish in $\mathbb{D}$, then for every real and  complex number $\alpha$ with $|\alpha|\leq1$ and $R\geq 1$,
         \begin{align*}
             |\tilde{\mathbb{S}}_a[P](Rz) - \alpha \tilde{\mathbb{S}}_a[P](z)| \leq & \frac{1}{2}\bigg[\biggl\{|R^n-\alpha||\tilde{\mathbb{S}}_a[z^n]|+ |1-\alpha|n|a|\biggr\}\max_{z \in B(\mathbb{D})}|P(z)|-\\& \biggl\{|R^n-\alpha||\tilde{\mathbb{S}}_a[z^n]|-|1-\alpha|n|a|\biggr\}\min_{z \in B(\mathbb{D})}|P(z)|\bigg].
         \end{align*}
         The above result is sharp and equality holds for $P(z) =z^n+1$.
     \end{corollary}
   \noindent  We derive the following result if we consider $\alpha = 0$ in inequality (\ref{2.3}).
     \begin{corollary}
          If $P(z)$ is a polynomial of degree $n$, which doesn't vanish in $\mathbb{D}$, then for every real and  complex number $ \beta $ with $|\beta|\leq 1$ and $R\geq 1$
     \begin{align*}
         \nonumber  \biggl|\Tilde{\mathbb{S}}_a[P]&(Rz)+\beta \biggl(\frac{R+1}{2}\biggr)^n \Tilde{\mathbb{S}}_a[P](z)\biggr| \leq \frac{1}{2}\biggl[ \biggl\{\biggl|R^n +\beta\bigg(\frac{R+1}{2}\bigg)^n\biggr||\Tilde{\mathbb{S}}_a[z^n]|+\\& \biggl|1 +\beta\biggl(\frac{R+1}{2}\biggr)^n \biggr|n|a|\biggr\}\max_{z \in B(\mathbb{D})}|P(z)|  -\biggl\{\bigg|R^n+\beta \bigg(\frac{R+1}{2}\bigg)^n\bigg||\Tilde{\mathbb{S}}_a[z^n]|- \\ &\bigg|1-\beta\bigg(\frac{R+1}{2}\bigg)^n\bigg|n|a|\biggr\}\min_{z\in B(\mathbb{D})}|P(z)|\bigg].
     \end{align*}
     The above result is sharp and holds for $P(z)=z^n+1$.
     \end{corollary}
   \noindent The next corollary is obtained by taking $\alpha=1$ in Theorem \ref{thm2}, dividing by $R-1$, and then letting $R\to 1$. This gives a refinement of Corollary \ref{cor2.1} for polynomials not vanishing in the unit disk.
     \begin{corollary}
            If $P(z)$ is a polynomial of degree $n$, which doesn't vanish in $\mathbb{D}$, then for every real and complex number $ \beta $ with $|\beta|\leq 1$ and $R\geq 1$
            \begin{align}{\label{2.4}}
                \nonumber  \bigg|z\tilde{\mathbb{S}}_a[P'](z)+\frac{n}{2}  \beta   \tilde{\mathbb{S}}_a[P](z) +P'(z)\biggr| \leq & \frac{n}{2}\biggl[ \biggl\{\bigg|1+\frac{\beta}{2}\bigg||\tilde{\mathbb{S}}_a[z^n]|+\frac{n}{2}|\beta||a|\biggr\} \max_{z \in B(\mathbb{D})}|P(z)|\\& - \biggl\{\bigg|1+\frac{\beta}{2}\bigg||\tilde{\mathbb{S}}_a[z^n]|-\frac{n}{2}|\beta||a|\biggr\}\min_{z \in B(\mathbb{D})}|P(z)|\biggr].
            \end{align}
            The result is the best possible and equality in \textnormal{(\ref{2.4})} holds for $P(z)=z^n+1$.
     \end{corollary}
     \begin{remark}
 If we consider $\alpha=\beta=0 $ and $ R=1$ in (\ref{2.3}), it reduces to the inequality (\ref{10}). Further, if we choose $a=-\frac{1}{z}$ in (\ref{2.3}), we get inequality (\ref{1.15}).
     \end{remark}
\begin{center}
\section*{LEMMA}	
\end{center}
 For the proof of these theorems, we need the following lemmas.
     The first lemma is due to Aziz \cite{aziz1987growth}.
\begin{lem}{\label{lem1}}
 If $P(z)$ is a polynomial of degree $n$, having all its zeros in $|z| \leq k, k\leq 1$ , then every $R\geq 1$ 
 \begin{align}{\label{3.1}}
     |P(Rz)| \geq \left(\frac{R+k}{1+k}\right)^n|P(z)|, \ for\ z \in B(\mathbb{D}). 
 \end{align}
\end{lem}
\begin{lem}{\label{lem2}}
 Let $P(z)$ be a polynomial of degree $n$ with all its zeros in $\overline{\mathbb{D}}$. If $a\in B(\mathbb{D})$ is not the exceptional value for $P$, then all the zeros of $\Tilde{\mathbb{S}}_a[P](z)$ lie in ${\overline{\mathbb{D}}}$.  
\end{lem}
\noindent The above lemma is due to Genenkova and Starkov \cite{ganenkova2019variations}. In addition, the next two lemmas are given by Deepak et al. \cite{kumar2025some}
   \begin{lem}{\label{lem3}}
   Let $P(z)$ and $F(z)$ be two polynomials such that $degP(z)\leq degF(z) =  n$. If  $F(z)$ has all zeros in $\overline{\mathbb{D}}$ and $|P(z)| \leq |F(z)|$, for $z\in B(\mathbb{D})$. Then for every real and complex number $\alpha, \beta $ with $|\alpha|\leq 1, |\beta| \leq 1$, $R\geq 1$ and $z \in \mathbb{C\backslash D}$
   \begin{align}{\label{3.2}}
      \nonumber \biggl|\Tilde{\mathbb{S}}_a&[P](Rz)-\alpha \Tilde{\mathbb{S}}_a[P](z)+\beta \biggl\{\biggl(\frac{R+1}{2}\biggr)^n -|\alpha|\biggr\}\Tilde{\mathbb{S}}_a[P](z)\biggl|\\ &\leq \biggl|\Tilde{\mathbb{S}}_a[F](Rz)-\alpha \Tilde{\mathbb{S}}_a[F](z)+\beta \biggl\{\biggl(\frac{R+1}{2}\biggr)^n -|\alpha|\biggr\}\Tilde{\mathbb{S}}_a[F](z)\biggr|.
   \end{align}
    \end{lem}
      \begin{lem}{\label{lem4}}
        If $P(z)$ is a polynomial of degree $n$, then for every real and complex number $\alpha, \beta$ with $|\alpha|\leq1,|\beta|\leq 1$ and $R\geq 1$,
     \begin{align}{\label{3.3}}
         \nonumber& \biggl|\Tilde{\mathbb{S}}_a[P](Rz)-\alpha\Tilde{\mathbb{S}}_a[P](z)+\beta\biggl\{\biggl(\frac{R+1}{2}\biggr)^n-|\alpha|\biggr\}\Tilde{\mathbb{S}}_a[P](z)\biggr|+\\& \nonumber \biggl|\Tilde{\mathbb{S}}_a[Q](Rz)-\alpha\Tilde{\mathbb{S}}_a[Q](z)+\beta\biggl\{\bigg(\frac{R+1}{2}\biggr)^n-|\alpha|\biggr\}\Tilde{\mathbb{S}}_a[Q](z)\biggr|  \\& \nonumber \leq \biggl[ \biggl|R^n-\alpha+\beta\biggl\{\bigg(\frac{R+1}{2}\bigg)^n-|\alpha|\biggr\}\biggr||\Tilde{\mathbb{S}}_a[z^n]|+ \\&\qquad \biggl|1-\alpha +\beta\bigg\{\biggl(\frac{R+1}{2}\biggr)^n-|\alpha|\bigg\}\biggr|n|a|\biggr]\max_{z\in B(\mathbb{D})}|P(z)|,\ for \ z \in \mathbb{C\backslash D}
        \end{align}
   where $Q(z)=z^n \overline{P(\frac{1}{\overline{z}})}$.
   \end{lem}
\begin{center}
\section*{PROOF OF THEOREMS}	
\end{center}
  \begin{proof}[\textbf{Proof of Theorem \ref{thm1}}]
     For $R =1$, the result is obvious. Therefore, we shall prove the result for $R>1$. If $P(z)$ has a zero in $B(\mathbb{D})$, then inequality (\ref{2.1}) is trivial. So we suppose that $P(z)$ has all its zeros in $\mathbb{D}$. If $m=\min_{z \in B(\mathbb{D})}|P(z)|$, then $0 < m \leq |P(z)|$ for $z \in B(\mathbb{D})$. Therefore, if $\lambda$ is a real and complex number such that $|\lambda| < 1$, then it follows by Rouche's theorem that the polynomial $F(z)=P(z)-\lambda mz^n$ of degree $n$, has all its zeros in $\mathbb{D}$.
     Applying Lemma \ref{lem1} to the polynomial $F(z)$ with $k=1$ and $R>1$, we get 
     \begin{align}{\label{4.1}}
         |F(Rz)| \geq \left(\frac{R+1}{2}\right)^n|F(z)| \ for \ z \in B(\mathbb{D}).
     \end{align}
     Hence for every real and complex number $\alpha$ with $|\alpha| \leq1$, we have 
     \begin{align}{\label{4.2}}
    \nonumber|F(Rz)-\alpha F(z)| & \geq |F(Rz)|-|\alpha||F(z)| \\
    & > \biggl\{\left(\frac{R+1}{2}\right)^n- |\alpha|\biggr\}|F(z)| \ for \ z \in B(\mathbb{D}).
     \end{align}
     Since $F(Re^{ i\theta})  \neq 0$ and from inequality (\ref{4.1}), we obtain 
     \[|F(e^{i \theta})| < |F(Re^{i \theta}|\]
     for every $R>1$ and $0\leq \theta \leq 2\pi$. Equivalently,
     \[|F(z)| <|F(Rz)| \ for \ z\in B(\mathbb{D}) \ and \ R>1.\]
     Since $F(Rz)$ has all its zeros in $|z| < \frac{1}{R} < 1$, a direct application of Rouche's theorem shows that for every real and complex number $\alpha$ with $|\alpha| \leq 1$, the polynomial $F(Rz)-\alpha F(z)$ has all its zeros in $\mathbb{D}$. Applying again Rouche's theorem, it follows from (\ref{4.2}) that for every real and complex number $\beta$ with $|\beta|\leq 1$ and $R > 1$ all the zeros of the polynomial
     \[T(z)=F(Rz)-\alpha F(z) +\beta \biggl\{\left(\frac{R+1}{2}\right)^n-|\alpha|\biggr\}F(z)\] lie in $\mathbb{D}$. So by Lemma \ref{lem2}, it follows that all the zeros of $\tilde{\mathbb{S}}_a[T](z)$ lie in $\mathbb{D}$. Replacing $F(z)$ by $P(z)-\lambda m z^n$, we conclude that for every real and complex number $\alpha, \beta$ and $\lambda$ with $|\lambda|< 1, |\alpha|\leq 1$ and $|\beta|\leq 1 $
     \begin{align}{\label{4.3}}
         \nonumber &\Tilde{\mathbb{S}}_a[P](Rz)- \alpha \Tilde{\mathbb{S}}_a[P](z)+\beta \biggl\{\biggl(\frac{R+1}{2}\biggr)^n-|\alpha| \biggr\}\Tilde{\mathbb{S}}_a[P](z) - \\& \lambda m \biggl[R^n-\alpha + \beta\biggl\{\biggl(\frac{R+1}{2}\biggr)^n -|\alpha|\biggr\}\biggr] \Tilde{\mathbb{S}}_a[z^n] \neq 0,\ for \  z \in \mathbb{C \backslash D}.
     \end{align} 
     This implies for every $\alpha, \beta$ with $|\alpha| \leq 1, |\beta|\leq 1$
     \begin{align}{\label{4.4}}
        \nonumber & \biggl|\Tilde{\mathbb{S}}_a[P](Rz)- \alpha \Tilde{\mathbb{S}}_a[P](z)+\beta \biggl\{\biggl(\frac{R+1}{2}\biggr)^n-|\alpha| \biggr\}\Tilde{\mathbb{S}}_a[P](z)\biggr| \\ &\geq m\biggl|R^n-\alpha + \beta\biggl\{\biggl(\frac{R+1}{2}\biggr)^n -|\alpha|\biggr\}\biggr| |\Tilde{\mathbb{S}}_a[z^n]|\ for \  z \in \mathbb{C \backslash D}. 
     \end{align}
     If it is not true, then there is a point $z=z_0 \in \mathbb{C \backslash D}$ such that for $R>1$
     \begin{align*}
    &\biggl|\tilde{\mathbb{S}}_a[P](Rz_0)-\alpha \tilde{\mathbb{S}}_a[P](z_0)+\beta\biggl\{\left(\frac{R+1}{2}\right)^n-|\alpha|\biggr\}\tilde{\mathbb{S}}_a[P](z_0)\biggr|\\& < m\biggl|R^n-\alpha + \beta\biggl\{\biggl(\frac{R+1}{2}\biggr)^n -|\alpha|\biggr\}\biggr| |\Tilde{\mathbb{S}}_a[z_0^n]|.
     \end{align*}
     We take
     \[\lambda = \frac{\tilde{\mathbb{S}}_a[P](Rz_0)-\alpha \tilde{\mathbb{S}}_a[P](z_0)+\beta\biggl\{\left(\frac{R+1}{2}\right)^n-|\alpha|\biggr\}\tilde{\mathbb{S}}_a[P](z_0)}{m\biggl[R^n-\alpha + \beta\biggl\{\biggl(\frac{R+1}{2}\biggr)^n -|\alpha|\biggr\}\biggr] \Tilde{\mathbb{S}}_a[z_0^n]},\] Then, $|\lambda| < 1$, and with this choice of $\lambda$, it follows from (\ref{4.3}) that $\tilde{\mathbb{S}}_a[T](z_0) = 0$ for $z = z_0 \in \mathbb{C\backslash D}$. However, this contradicts the fact that all zeros of $\tilde{\mathbb{S}}_a[T](z)$ lie in $\mathbb{D}$. Consequently, in particular, (\ref{4.4}) holds for $\alpha, \beta$ satisfying $|\alpha| \leq 1, |\beta| \leq 1$, and $R \geq 1$.
     \begin{align*}
       &\biggl| \Tilde{\mathbb{S}}_a[P](Rz)-\alpha \Tilde{\mathbb{S}}_a[P](z)+\beta \biggl\{\biggl(\frac{R+1}{2}\biggr)^n-|\alpha| \biggr\}\Tilde{\mathbb{S}}_a[P](z)\biggr|\\& \geq \biggl|R^n-\alpha +\beta\biggl\{\biggl(\frac{R+1}{2}\biggr)^n -|\alpha|\biggr\}\biggr||\tilde{\mathbb{S}}_a[z^n]|\min_{z\in B(\mathbb{D})}|P(z)|.
     \end{align*}
     This completes the proof of Theorem \ref{thm1}.
     \end{proof}
 \begin{proof}[\textbf{Proof of Theorem \ref{thm2}}]
      For $R = 1$, there is nothing to prove. Therefore, we assume that $R>1$. By hypothesis the polynomial $P(z)\neq 0$ in $\mathbb{D}$, if $m= \min_{z \in B(\mathbb{D})} |P(z)|$, then $m \leq |P(z)|$ for $z \in B(\mathbb{D})$. Therefore, for every real and complex number $\gamma$ with $|\gamma|< 1$, it follows by Rouche's theorem that the polynomial $G(z)=P(z)-\gamma m$ has no zeros in $\mathbb{D}$. Now if \[H(z) = z^n \overline{G\bigg(\frac{1}{\overline{z}}\bigg)}=Q(z)-m\overline{\gamma}z^n,\] then all the zeros of $H(z)$ lie in $\mathbb{D}$ and $|G(z)|= |H(z)|$ for $z \in B(\mathbb{D})$. 
      Therefore by lemma \ref{lem3}, we have for $|\alpha|\leq 1,|\beta|\leq 1$ and $z \in B(\mathbb{D})$
\begin{align}{\label{4.5}}
   \nonumber& \biggl|\Tilde{\mathbb{S}}_a[P](Rz)-\alpha\Tilde{\mathbb{S}}_a[P](z)+\beta\biggl\{\bigg(\frac{R+1}{2}\biggr)^n-|\alpha|\biggr\}\Tilde{\mathbb{S}}_a[P](z)+\\ & \nonumber\quad\gamma\biggr[1-\alpha +\beta\bigg\{\biggl(\frac{R+1}{2}\biggr)^n-|\alpha|\bigg\}\bigg]nam\biggr| \\ &\nonumber \leq\biggl|\Tilde{\mathbb{S}}_a[Q](Rz)-\alpha\Tilde{\mathbb{S}}_a[Q](z)+\beta\biggl\{\bigg(\frac{R+1}{2}\biggr)^n-|\alpha|\biggr\}\Tilde{\mathbb{S}}_a[Q](z)-\\ & \quad \overline{\gamma}\bigg[R^n-\alpha+\beta\biggl\{\bigg(\frac{R+1}{2}\bigg)^n-|\alpha|\biggr\}\biggr]m\Tilde{\mathbb{S}}_a[z^n]\biggr|,
\end{align}
 where  $Q(z)=z^n\overline{P(\frac{1}{\overline{z}})}$.\\
Since all the zeros of $Q(z)$ lie in $\mathbb{D}$ and $m =\min_{z \in B(\mathbb{D})}|P(z)|= \min_{z \in B(\mathbb{D})}|Q(z)|$, by applying Theorem \ref{thm1}, we obtain 
\begin{align*}
    &\biggl|\Tilde{\mathbb{S}}_a[Q](Rz)-\alpha\Tilde{\mathbb{S}}_a[Q](z)+\beta\biggl\{\bigg(\frac{R+1}{2}\biggr)^n-|\alpha|\biggr\}\Tilde{\mathbb{S}}_a[Q](z)\biggr|\\& \geq \biggl|R^n-\alpha +\beta \biggl\{\bigg(\frac{R+1}{2}\bigg)^n-|\alpha|\biggr\}\biggr||\tilde{\mathbb{S}}_a[z^n]|\min_{z \in B(\mathbb{D})}|Q(z)|
\\& = \biggl|R^n-\alpha +\beta \biggl\{\bigg(\frac{R+1}{2}\bigg)^n-|\alpha|\biggr\}\biggr|m |\tilde{\mathbb{S}}_a[z^n]|,
\end{align*}
for $|\alpha|\leq 1, |\beta| \leq 1$.\\
Now choosing the argument of $\gamma$ in (\ref{4.5}) on the right hand, we get
\begin{align}{\label{4.6}}
    \nonumber&\biggl|\Tilde{\mathbb{S}}_a[Q](Rz)-\alpha\Tilde{\mathbb{S}}_a[Q](z)+\beta\biggl\{\bigg(\frac{R+1}{2}\biggr)^n-|\alpha|\biggr\}\Tilde{\mathbb{S}}_a[Q](z)-\\& \nonumber\overline{\gamma}\bigg[R^n-\alpha+\beta\biggl\{\bigg(\frac{R+1}{2}\bigg)^n-|\alpha|\biggr\}m\Tilde{\mathbb{S}}_a[z^n]\biggr]\biggr| \\& \nonumber = \biggl|\Tilde{\mathbb{S}}_a[Q](Rz)-\alpha\Tilde{\mathbb{S}}_a[Q](z)+\beta\biggl\{\bigg(\frac{R+1}{2}\biggr)^n-|\alpha|\biggr\}\Tilde{\mathbb{S}}_a[Q](z)\biggr|- \\&  \qquad m |\overline{\gamma}|\biggl|R^n-\alpha+\beta\biggl\{\bigg(\frac{R+1}{2}\bigg)^n-|\alpha|\biggr\}\biggr||\Tilde{\mathbb{S}}_a[z^n]|.
\end{align}
From inequalities (\ref{4.5}) and (\ref{4.6}) and taking $|\gamma| \to 1$, we have
 \begin{align*}
 &\biggl|\Tilde{\mathbb{S}}_a[P](Rz)-\alpha\Tilde{\mathbb{S}}_a[P](z)+\beta\biggl\{\bigg(\frac{R+1}{2}\biggr)^n-|\alpha|\biggr\}\Tilde{\mathbb{S}}_a[P](z)\biggr| -\\ & \quad m\bigg|1-\alpha+\beta\biggl\{\bigg(\frac{R+1}{2}\bigg)^n-|\alpha|\biggr\}\biggr|n|a|\\& \leq \biggl|\Tilde{\mathbb{S}}_a[Q](Rz)-\alpha\Tilde{\mathbb{S}}_a[Q](z)+\beta\biggl\{\bigg(\frac{R+1}{2}\biggr)^n-|\alpha|\biggr\}\Tilde{\mathbb{S}}_a[Q](z)\biggr|-\\ & \quad m\biggl|R^n-\alpha+\beta\biggl\{\bigg(\frac{R+1}{2}\bigg)^n-|\alpha|\biggr\}\biggr||\Tilde{\mathbb{S}}_a[z^n]|.
\end{align*}
Equivalently
\begin{align*}
     &\biggl|\Tilde{\mathbb{S}}_a[P](Rz)-\alpha\Tilde{\mathbb{S}}_a[P](z)+\beta\biggl\{\bigg(\frac{R+1}{2}\biggr)^n-|\alpha|\biggr\}\Tilde{\mathbb{S}}_a[P](z)\biggr|\\ & \leq \biggl|\Tilde{\mathbb{S}}_a[Q](Rz)-\alpha\Tilde{\mathbb{S}}_a[Q](z)+\beta\biggl\{\bigg(\frac{R+1}{2}\biggr)^n-|\alpha|\biggr\}\Tilde{\mathbb{S}}_a[Q](z)\biggr|-\\ & \quad m\bigg[\biggl|R^n-\alpha+\beta\biggl\{\bigg(\frac{R+1}{2}\bigg)^n-|\alpha|\biggr\}\biggr||\Tilde{\mathbb{S}}_a[z^n]|- \\ & \qquad \biggl|1-\alpha+\beta\biggl\{\bigg(\frac{R+1}{2}\bigg)^n-|\alpha|\biggr\}\biggr|n|a|\biggr].
\end{align*}
Which implies for every real or complex number $\alpha, \beta$ with $ |\alpha| \leq 1, |\beta|\leq 1, R>1$ and $z \in B(\mathbb{D})$  
     \begin{align*}
         2\biggl|&\Tilde{\mathbb{S}}_a[P](Rz)-\alpha \Tilde{\mathbb{S}}_a[P](z)+\beta \biggl\{\biggl(\frac{R+1}{2}\biggr)^n -|\alpha|\biggr\}\Tilde{\mathbb{S}}_a[P](z)\biggr|\\ &\leq \biggl|\Tilde{\mathbb{S}}_a[Q](Rz)-\alpha \Tilde{\mathbb{S}}_a[Q](z)+\beta \biggl\{\biggl(\frac{R+1}{2}\biggr)^n -|\alpha|\biggr\}\Tilde{\mathbb{S}}_a[Q](z)\biggr|+\\& \qquad \biggl|\Tilde{\mathbb{S}}_a[P](Rz)-\alpha \Tilde{\mathbb{S}}_a[P](z)+ \beta \biggl\{\biggl(\frac{R+1}{2}\biggr)^n -|\alpha|\biggr\}\Tilde{\mathbb{S}}_a[P](z)\biggr| -\\ & \quad m\bigg[\biggl|R^n-\alpha+\beta\biggl\{\bigg(\frac{R+1}{2}\bigg)^n-|\alpha|\biggr\}\biggr||\Tilde{\mathbb{S}}_a[z^n]|- \\ & \qquad \biggl|1-\alpha+\beta\biggl\{\bigg(\frac{R+1}{2}\bigg)^n-|\alpha|\biggr\}\biggr|n|a|\biggr].
     \end{align*}
     Applying Lemma \ref{lem4}, we have
     \begin{align*}
        2\biggl|&\Tilde{\mathbb{S}}_a[P](Rz)-\alpha \Tilde{\mathbb{S}}_a[P](z)+\beta \biggl\{\biggl(\frac{R+1}{2}\biggr)^n -|\alpha|\biggr\}\Tilde{\mathbb{S}}_a[P](z)\biggr|\\& \leq\biggl[ \biggl|R^n-\alpha+\beta\biggl\{\bigg(\frac{R+1}{2}\bigg)^n-|\alpha|\biggr\}\biggr||\Tilde{\mathbb{S}}_a[z^n]|+\\& \qquad\biggl|1-\alpha +\beta\bigg\{\biggl(\frac{R+1}{2}\biggr)^n-|\alpha|\bigg\}\biggr|n|a|\biggr]\max_{z \in \mathbb{D}}|P(z)| - \\ & \quad m\bigg[\biggl|R^n-\alpha+\beta\biggl\{\bigg(\frac{R+1}{2}\bigg)^n-|\alpha|\biggr\}\biggr||\Tilde{\mathbb{S}}_a[z^n]|- \\ & \qquad \biggl|1-\alpha+\beta\biggl\{\bigg(\frac{R+1}{2}\bigg)^n-|\alpha|\biggr\}\biggr|n|a|\biggr].
     \end{align*}
     This is the complete proof of Theorem \ref{thm2}. 
 \end{proof}
 \bibliographystyle{amsplain}
 \bibliography{Derivative_Smirnov-II}
\end{document}